\documentclass[reqno]{amsart}
\usepackage{amssymb}
\usepackage{epsf}

\newtheorem{theorem}{Theorem}

\theoremstyle{definition}
\newtheorem*{definition}{Definition}

\theoremstyle{remark}

\begin{document}

\title[A Bailey Tree for Integrals]
{A Bailey tree for integrals}

\author{V.P. Spiridonov}

\address{Bogoliubov Laboratory of Theoretical Physics, JINR,
Dubna, Moscow reg. 141980, Russia}

\thanks{{\em Date:} May 2003; to be published in {\em Theor. Math.
Phys.} (Russian)}

\begin{abstract}
The notion of integral Bailey pairs is introduced.
Using the single variable elliptic beta integral, we
construct an infinite binary  tree of identities for
elliptic hypergeometric integrals. Two particular
sequences of identities are explicitly described.
\end{abstract}

\maketitle

Series and integrals of hypergeometric type have many
applications in mathematical physics. Investigation
of the corresponding special functions can therefore be relevant
for various computations in theoretical models of reality.
Bailey chains provide powerful tools for generating infinite
sequences of summation or transformation formulas for
hypergeometric type series. For a review of the corresponding formalism
and many $q$-series identities obtained with its help, see, e.g.,
\cite{and:bailey}, \cite{war:50}. Recently, we extended
the scope of applications of this technique to elliptic
hypergeometric series \cite{spi:bailey}. More
precisely, an elliptic generalization of the Andrews' well-poised
Bailey chain \cite{and:bailey} has been found and its simplest
consequences have been analyzed. This lifted also one
of the Andrews-Berkovich identities \cite{and-ber:bailey}
to the elliptic level.

During the work on \cite{spi:bailey}, we came to a
principal conclusion that there must exist Bailey chains for
integrals, but the first attempts to build a simple example
of such a chain were not successful. In \cite{spi:integrals},
a Bailey type symmetry transformation was constructed for
a pair of elliptic hypergeometric integrals. This result gives
the tools necessary for an appropriate generalization
of the elliptic Bailey chain technique of \cite{spi:bailey}.
In the present note, we construct two examples
of Bailey chains for integrals directly at the level of the
most complicated known type of beta integrals of one variable,
namely, the elliptic beta integral of \cite{spi:elliptic}.
Being used together, they form a binary tree of identities
for elliptic hypergeometric integrals.

We denote by $\mathbb{T}$ the positively oriented unit circle
and take two base variables $q, p$, $|q|,|p|<1$, and five complex
parameters $t_m, m=0, \dots,4,$ satisfying the inequalities
$|t_m|<1,\; |pq|<|A|,$ where $A\equiv\prod_{r=0}^4t_r.$
Then the elliptic beta integral
\begin{equation}\label{ell-int}
\mathcal{N}_E(\mathbf{t})=
\frac{1}{2\pi i}\int_\mathbb{T}\frac{\prod_{m=0}^4\Gamma(z^\pm t_m)}
{\Gamma(z^{\pm 2},z^\pm A)}\frac{dz}{z},
\end{equation}
admits the exact evaluation \cite{spi:elliptic}:
\begin{equation}\label{result}
\mathcal{N}_E(\mathbf{t})=\frac{2\prod_{0\leq m<s\leq 4}
\Gamma(t_mt_s)}
{(q;q)_\infty(p;p)_\infty\prod_{m=0}^4\Gamma(At_m^{-1})}.
\end{equation}
The standard notation $(q;q)_\infty=\prod_{n=1}^\infty(1-q^n)$ and
the following conventions are implied throughout the paper:
\begin{eqnarray*}
&& \Gamma(tz^\pm)=\Gamma(tz,tz^{-1}),\quad
\Gamma(z^{\pm 2})=\Gamma(z^2,z^{-2}),
\\
&& \Gamma(tz^\pm x^\pm)=\Gamma(tzx,tzx^{-1},tz^{-1}x,tz^{-1}x^{-1}),
\\
&& \Gamma(t_1,\ldots,t_m)\equiv\prod_{r=1}^m \Gamma(t_r;q,p),
\end{eqnarray*}
where $\Gamma(z;q,p)$ is the elliptic gamma function
\begin{equation}
\Gamma(z;q,p) = \prod_{j,k=0}^\infty\frac{1-z^{-1}q^{j+1}p^{k+1}}
{1-zq^jp^k}
\label{ell-gamma}\end{equation}
(for a brief historical account of this function and relevant references,
see, e.g., \cite{spi:integrals}).
For $p=0$, formula (\ref{ell-int}) is reduced to the Nassrallah-Rahman
$q$-beta integral which can be reduced further to the
Askey-Wilson integral \cite{gas-rah:basic}.

The main definition of the present paper introduces the notion
of Bailey pairs for integrals.

\begin{definition}
Two functions $\alpha(z,t)$ and $\beta(z,t)$, $z, t\in \mathbb{C}$,
are said to form {\em an integral elliptic Bailey pair} with respect
to the parameter $t$ if they satisfy the following relation
\begin{equation}
\beta(w,t)=\kappa
\int_{\mathbb{T}}\Gamma(tw^\pm z^\pm)\, \alpha(z,t)\frac{dz}{z},
\qquad \kappa\equiv\frac{(p;p)_\infty (q;q)_\infty}{4\pi i}.
\label{bint-pair}\end{equation}
\end{definition}

Suppose that we have found a particular integral elliptic Bailey pair.
Then there exists a possibility to build an infinite
sequence of such pairs out of the given one.

\begin{theorem}
Whenever two functions $\alpha(z,t)$ and $\beta(z,t)$
form an integral elliptic Bailey pair with respect to $t$,
the new functions defined for $w\in\mathbb{T}$ as
\begin{equation}
\alpha'(w,st)=\frac{\Gamma(tuw^\pm)}
{\Gamma(ts^2uw^\pm)}\,\alpha(w,t)
\label{a'}\end{equation}
and
\begin{equation}
\beta'(w,st)=\kappa
\frac{\Gamma(t^2s^2,t^2suw^\pm)}{\Gamma(s^2,t^2,suw^\pm)}
\int_{\mathbb{T}} \frac{\Gamma(sw^\pm x^\pm,ux^\pm)}
{\Gamma(x^{\pm 2},t^2s^2ux^\pm)}\,\beta(x,t)\frac{dx}{x}
\label{b'}\end{equation} 
form an integral elliptic Bailey pair with respect
to the parameter $st$. Here $t,s,u$ are arbitrary complex parameters
satisfying the constraints $|t|,|s|,|u|<1, |pq|<|t^2s^2u|$.
\end{theorem}

\begin{proof}
We substitute definition (\ref{bint-pair}) into relation (\ref{b'})
\begin{eqnarray}\nonumber
\lefteqn{\beta'(w,st)=\kappa^2
\frac{\Gamma(t^2s^2,t^2suw^\pm)}{\Gamma(s^2,t^2,suw^\pm)} } &&
\\ && \makebox[2em]{}
\times \int_{\mathbb{T}^2}
\frac{\Gamma(sw^\pm x^\pm,tz^\pm x^\pm,ux^\pm)}
{\Gamma(x^{\pm 2},t^2s^2ux^\pm)}\,\alpha(z,t)\frac{dz}{z}\frac{dx}{x}.
\label{proof}\end{eqnarray}
Since the integrands are bounded on $\mathbb{T}$, we may change in this
expression the order of integrations (i.e., integrate first
with respect to $x$) and apply formula (\ref{ell-int}). This yields
\begin{eqnarray}\nonumber
&& \beta'(w,st)=\kappa
\int_{\mathbb{T}}\frac{\Gamma(stw^\pm z^\pm,tuz^\pm)}
{\Gamma(ts^2uz^\pm)}\; \alpha(z,t)\frac{dz}{z}
\\ && \makebox[4em]{}
= \kappa
\int_{\mathbb{T}}\Gamma(stw^\pm z^\pm)\; \alpha'(z,st)\frac{dz}{z},
\label{st}\end{eqnarray}
that is the functions $\alpha'(w,st)$ and $\beta'(w,st)$ form an integral
Bailey pair with respect to the parameter $st$.
\end{proof}
In some sense, this theorem may be called as an integral Bailey lemma.

If we substitute into integral (\ref{ell-int}) relations
$t_3=tw, t_4=tw^{-1}$,
then it is not difficult to figure out an initial integral elliptic
Bailey pair:
\begin{eqnarray}\label{a_1}
&& \alpha(z,t)=\frac{\prod_{r=0}^2\Gamma(t_rz^\pm)}
{\Gamma(z^{\pm 2}, t^2t_0t_1t_2z^\pm)},
\\ &&
\beta(w,t)=\Gamma(t^2)\prod_{0\leq r<j\leq 2}\frac{\Gamma(t_rt_j)}
{\Gamma(t^2t_rt_j)}\,
\frac{\prod_{r=0}^2\Gamma(tt_rw^\pm)}
{\Gamma(tt_0t_1t_2w^\pm)},
\label{b_1}\end{eqnarray}
where $|t|, |t_r|<1, |pq|<|t^2t_0t_1t_2|$.
Substituting these expressions into chain rules (\ref{a'}) and
(\ref{b'}), we obtain
\begin{eqnarray*}
&& \alpha'(z,st)=\frac{\Gamma(tuz^\pm)\prod_{r=0}^2\Gamma(t_rz^\pm)}
{\Gamma(ts^2uz^\pm,z^{\pm2},t^2t_0t_1t_2z^\pm)},
\\
&& \beta'(w,st)=\kappa \frac{\Gamma(t^2s^2,t^2suw^\pm)}
{\Gamma(s^2,suw^\pm)}
\prod_{0\leq r<j\leq 2}\frac{\Gamma(t_rt_j)}{\Gamma(t^2t_rt_j)}
\\
&& \makebox[4em]{} \times
\int_{\mathbb{T}}
\frac{\Gamma(sw^\pm x^\pm,ux^\pm)\prod_{r=0}^2\Gamma(tt_rx^\pm)}
{\Gamma(x^{\pm 2},t^2s^2ux^\pm,tt_0t_1t_2x^\pm)}\frac{dx}{x}.
\end{eqnarray*}
Being plugged into key generating relation (\ref{st}), these
functions yield a symmetry transformation for
two elliptic hypergeometric integrals:
\begin{eqnarray}\nonumber
&& \prod_{j=0}^{2}\frac{\Gamma(Bt_j^{-1})}{\Gamma(t^2Bt_j^{-1})}
\int_{\mathbb{T}}\frac{\prod_{j=1}^{3}
\Gamma(tt_jz^\pm,s_jz^\pm)}
{\Gamma(z^{\pm2},t^2Sz^\pm,tBz^\pm)}\frac{dz}{z}
\\ && \makebox[2em]{}
=\prod_{j=0}^{2}\frac{\Gamma(Ss_j^{-1})}{\Gamma(t^2Ss_j^{-1})}
\int_{\mathbb{T}}\frac{\prod_{j=1}^{3}
\Gamma(ts_jz^\pm,t_jz^\pm)}
{\Gamma(z^{\pm2},t^2Bz^\pm,tSz^\pm)}\frac{dz}{z},
\label{an-trans}\end{eqnarray}
where $s_0=u,s_1=sw,s_2=sw^{-1}$ and
$B=\prod_{j=0}^{2}t_j, S=\prod_{j=0}^{2}s_j$.
Identity (\ref{an-trans}) was established in \cite{spi:integrals}
and the integral elliptic Bailey chain (\ref{bint-pair})-(\ref{b'})
is constructed with the help of a generalization of the corresponding method of
derivation of this equality.

Iterative application of rules (\ref{a'}) and (\ref{b'})
generates an infinite chain of identities for elliptic hypergeometric
integrals. The result of $m$-fold iteration has the form
\begin{eqnarray}\label{a-m}
&& \alpha^{(m)}(x,\prod_{l=1}^ms_l\, t)=
\prod_{k=1}^m\frac{\Gamma(t\prod_{l=1}^{k-1}s_l\, u_kx^\pm)}
{\Gamma(t\prod_{l=1}^{k-1}s_l\, s_k^2u_kx^\pm)} \alpha^{(0)}(x,t),
\\ \label{b-m}  &&
\beta^{(m)}(x_{m+1},\prod_{l=1}^ms_l\, t)=\kappa^m
\prod_{k=1}^m\frac{\Gamma(t^2\prod_{l=1}^{k}s_l^2)}
{\Gamma(s_k^2,t^2\prod_{l=1}^{k-1}s_l^2)}
\\ && \makebox[2em]{} \times
\int_{\mathbb{T}^m} \prod_{k=1}^m\frac
{\Gamma(t^2\prod_{l=1}^{k-1}s_l^2\, s_ku_kx_{k+1}^\pm,
s_kx_{k+1}^\pm x_k^\pm,u_kx_k^\pm)}
{\Gamma(s_ku_kx_{k+1}^\pm,x_k^{\pm 2}, t^2\prod_{l=1}^ks_l^2\,
u_kx_k^\pm)} \beta^{(0)}(x_1,t)\prod_{l=1}^m\frac{dx_l}{x_l}.
\nonumber\end{eqnarray}

Using (\ref{a_1}) and (\ref{b_1}) as the initial
functions  $\alpha^{(0)}(x,t)$ and $\beta^{(0)}(x,t)$
and substituting (\ref{a-m}) and (\ref{b-m}) into (\ref{bint-pair}),
we obtain the identity
\begin{eqnarray}\nonumber
\lefteqn{ \kappa^{m-1}\prod_{k=1}^m
\frac{\Gamma(t^2\prod_{l=1}^{k}s_l^2)}
{\Gamma(s_k^2,t^2\prod_{l=1}^{k-1}s_l^2)}
\int_{\mathbb{T}^m}\frac{\prod_{r=0}^2\Gamma(tt_rx_1^\pm)}
{\Gamma(tt_0t_1t_2x_1^\pm)}  } &&
\\ \nonumber && \makebox[2em]{} \times
\prod_{k=1}^m\frac
{\Gamma(t^2\prod_{l=1}^{k-1}s_l^2\, s_ku_kx_{k+1}^\pm,
s_kx_{k+1}^\pm x_k^\pm,u_kx_k^\pm)}
{\Gamma(s_ku_kx_{k+1}^\pm,x_k^{\pm 2}, t^2\prod_{l=1}^ks_l^2\,
u_kx_k^\pm)} \prod_{l=1}^m\frac{dx_l}{x_l}
\\ \nonumber &&
=\frac{1}{\Gamma(t^2)}\prod_{0\leq r<j\leq 2}\frac{\Gamma(t^2t_rt_j)}
{\Gamma(t_rt_j)} \int_{\mathbb{T}}
\frac{\Gamma(t\prod_{l=1}^ms_l\, x_{m+1}^\pm x^\pm)
\prod_{r=0}^2\Gamma(t_rx^\pm)}{\Gamma(x^{\pm 2},t^2t_0t_1t_2x^\pm)}
\\ && \makebox[2em]{} \times
\prod_{k=1}^m \frac{\Gamma(t\prod_{l=1}^{k-1}s_l\, u_kx^\pm)}
{\Gamma(t\prod_{l=1}^{k-1}s_l\, s_k^2u_kx^\pm)}\frac{dx}{x},
\label{id-seq}\end{eqnarray}
which can be rewritten as (\ref{an-trans}) for $m=1$.
In some sense, relation (\ref{id-seq}) is a high level
integral generalization of the multiple series Rogers-Ramanujan
type identities derived in \cite{and:multiple}. Integral analogues of
some of the Rogers-Ramanujan identities have been considered
in \cite{gis:variants}. Taking into account a great deal of activity
around Bailey chains for series, it is surprising that their integral
analogues did not get appropriate attention in the literature.
Some examples of the one step Bailey type transformations for
integrals were known for a long time. Actually, all the changes of
the orders of integrations in multiple integrals used in
\cite{and:short}, \cite{den-gus:beta}, \cite{gus:some2},
\cite{gus-rak:beta}, \cite{die-spi:selberg}, \cite{spi:integrals}
may be interpreted as integral generalizations of the
Bailey's change of the order of summations for series
\cite{and:bailey}, \cite{war:50}.
The general idea of using a change of the integrations order for obtaining
useful information about integrals goes back to Poisson's and Jacobi's
methods of derivation of the ordinary beta integral \cite{aar:special}.

It is known that for a taken definition of Bailey pairs there sometimes
exist several different chain substitution rules
(or Bailey lemmas). In this case one gets a more complicated
structure of the derived set of identities, which is called as the
Bailey lattice \cite{aab:bailey} or tree \cite{and:bailey}. It appears
that in our case we can extend the integral Bailey chain described in
Theorem 1 to a tree due to a complimentary integral
Bailey lemma described below.

\begin{theorem}
Whenever two functions $\alpha(z,t)$ and $\beta(z,t)$
form an integral elliptic Bailey pair with respect to $t$,
$|t|<1$, the new functions defined for $w\in\mathbb{T}$ as
\begin{eqnarray}
\alpha'(w,t)=\kappa \frac{\Gamma(s^2t^2,uw^\pm)}
{\Gamma(s^2,t^2,w^{\pm2},t^2s^2uw^{\pm})}
\int_{\mathbb{T}} \frac{\Gamma(t^2sux^\pm,sw^\pm x^\pm)}
{\Gamma(sux^\pm)}\,\alpha(x,st)\frac{dx}{x}
\label{a2}\end{eqnarray}
and
\begin{equation}
\beta'(w,t)=\frac{\Gamma(tuw^\pm)}
{\Gamma(ts^2uw^\pm)}\,\beta(w,st)
\label{b2}\end{equation}
form an integral elliptic Bailey pair with respect
to $t$ as well. Here $s$ and $u$ are arbitrary complex parameters (they are
not related to $s,u$ appearing in Theorem 1) satisfying
the constraints $|s|,|u|<1, |pq|<|t^2s^2u|$.
\end{theorem}

\begin{proof}
Verification of this statement follows the proof of the previous theorem.
For this it is necessary to replace $\alpha(w,t)$ and $\beta(w,t)$ in
equality (\ref{bint-pair}) by primed expressions, to substitute relation
(\ref{a2}) into it, to change the order of integrations, and to apply
formula (\ref{result}). We skip the details of consideration for their
simplicity.
\end{proof}

In \cite{spi:bailey}, only one Bailey lemma was established for
elliptic hypergeometric series. Our Theorem 1 was derived in a
heuristic analogy with that. It is natural to expect that there
exists a limiting process which would allow one to pass from Theorem 1
to the corresponding Bailey chain for series. Evidently, our
Theorem 2 should have a series analogue as well which would lead to
an elliptic Bailey tree for series (its analogue for $q$-hypergeometric
series would be qualitatively different from
the tree considered in \cite{and:bailey}, \cite{and-ber:bailey}).

As an illustration, we describe one particular sequence of identities
generated by the derived Bailey tree. For this we take functions
(\ref{a_1}) and (\ref{b_1}) and apply to them
transformations (\ref{a2}) and (\ref{b2}) with parameters $s=s_1,u=u_1$.
Then we apply to the resulting $\alpha'(w,t)$ and $\beta'(w,t)$
transformations (\ref{a'}) and (\ref{b'}) with the parameters $s=s_2,u=u_2$.
As a result, we obtain
\begin{eqnarray}\nonumber
&& \alpha''(w,s_2t)=\kappa
\frac{\Gamma(t^2s_1^2,tu_2w^\pm,u_1w^\pm)}
{\Gamma(s_1^2,t^2,w^{\pm2},ts_2^2u_2w^\pm,t^2s_1^2u_1w^\pm)}
\\ && \makebox[5em]{} \times
\int_{\mathbb{T}} \frac{\Gamma(t^2s_1u_1x^\pm,s_1w^\pm x^\pm)
\prod_{r=0}^2\Gamma(t_rx^\pm)}{\Gamma(x^{\pm2},s_1u_1x^\pm,
t^2s_1^2t_0t_1t_2x^\pm)}\frac{dx}{x},
\\ \nonumber &&
\beta''(w,s_2t)=\kappa
\frac{\Gamma(t^2s_2^2,t^2s_1^2,t^2s_2u_2w^\pm)}
{\Gamma(s_2^2,t^2,s_2u_2w^\pm)}\prod_{0\leq r<j\leq 2}
\frac{\Gamma(t_rt_j)}{\Gamma(t^2s_1^2t_rt_j)}
\\ && \makebox[5em]{} \times
\int_{\mathbb{T}}\frac{\Gamma(s_2w^\pm x^\pm,tu_1x^\pm,u_2x^\pm)
\prod_{r=0}^2\Gamma(ts_1t_r x^\pm)}{\Gamma(x^{\pm2},t^2s_2^2u_2x^\pm,
ts_1^2u_1x^\pm,ts_1t_0t_1t_2x^\pm)}\frac{dx}{x}.
\label{ab''}\end{eqnarray}
Substituting these functions into basic relation (\ref{bint-pair}) with $t$
replaced by $s_2t$, we obtain an equality for two integrals
\begin{eqnarray} \nonumber
&& \int_{\mathbb{T}}\frac{\Gamma(s_2w^\pm x^\pm,tu_1x^\pm,u_2x^\pm)
\prod_{r=0}^2\Gamma(ts_1t_r x^\pm)}{\Gamma(x^{\pm2},t^2s_2^2u_2x^\pm,
ts_1^2u_1x^\pm,ts_1t_0t_1t_2x^\pm)}\frac{dx}{x}
\\ && \nonumber \makebox[4em]{}
=\kappa \prod_{0\leq r<j\leq 2}\frac{\Gamma(t^2s_1^2t_rt_j)}{\Gamma(t_rt_j)}
\frac{\Gamma(s_2^2,s_2u_2w^\pm)}{\Gamma(s_1^2,t^2s_2^2,t^2s_2u_2w^\pm)}
\\ \nonumber && \makebox[5em]{}\times
\int_{\mathbb{T}^2}\Gamma(s_1x_2^\pm x_1^\pm)
\frac{\Gamma(ts_2w^\pm x_2^\pm,tu_2x_2^\pm,u_1x_2^\pm)}
{\Gamma(x_2^{\pm2},ts_2^2u_2x_2^\pm,t^2s_1^2u_1x_2^\pm)}
\\ && \makebox[5em]{}\times
\frac{\Gamma(t^2s_1u_1x_1^\pm)\prod_{r=0}^2
\Gamma(t_rx_1^\pm)}
{\Gamma(x_1^{\pm2},s_1u_1x_1^\pm,t^2s_1^2t_0t_1t_2x_1^\pm)}
\frac{dx_1}{x_1}\frac{dx_2}{x_2}.
\label{ident1}\end{eqnarray}
This identity was derived under the constraints that all parameters
lie inside the unit circle and $|t^2s_1^2t_0t_1t_2|, |t^2s_1^2u_2|, |t^2s_2^2u_2|
>|pq|$. However, by analyticity it remains true for a larger region of
values of parameters provided we replace $\mathbb{T}$ by contours of
integration which encircle sequences of poles of the integrands
converging to zero and exclude other ones.

It is not difficult to find the result of $m$-fold iteration of
transformations (\ref{a2}) and (\ref{b2}) with different parameters $u_k,s_k$:
\begin{eqnarray}\nonumber
&& \alpha^{(m)}(x_{m+1},t)=\kappa^m\prod_{k=1}^m
\frac{\Gamma(\prod_{l=k}^ms_l^2t^2)}
{\Gamma(s_k^2,\prod_{l=k+1}^ms_l^2t^2)}
\\  && \makebox[4em]{} \times
\int_{\mathbb{T}^m}\prod_{k=1}^m\frac{\Gamma(u_kx_{k+1}^\pm,
\prod_{l=k+1}^ms_l^2t^2s_ku_kx_k^\pm,s_kx_{k+1}^\pm x_k^\pm)}
{\Gamma(x_{k+1}^{\pm2},\prod_{l=k}^ms_l^2t^2u_kx_{k+1}^\pm,s_ku_kx_k^\pm)}
\nonumber \\  && \makebox[6em]{} \times
\alpha^{(0)}\left(x_1,\prod_{l=1}^ms_l\, t\right)
\frac{dx_1}{x_1}\cdots\frac{dx_m}{x_m},
\label{a_m}\\
&& \beta^{(m)}(w,t)=\prod_{k=1}^m\frac{\Gamma(t\prod_{l=k+1}^ms_lu_kw^\pm)}
{\Gamma(t\prod_{l=k+1}^ms_ls_k^2u_kw^\pm)}\,
\beta^{(0)}\left(w,\prod_{l=1}^ms_l\, t\right).
\label{b_m}\end{eqnarray}
An application of transformations (\ref{a'}) and (\ref{b'}) with the
parameters $s=s_{m+1}, u=u_{m+1}$ to these functions yields:
\begin{eqnarray}\nonumber
&& \alpha'(x_{m+1},s_{m+1}t)=\frac{\Gamma(tu_{m+1}x_{m+1}^\pm)}
{\Gamma(ts_{m+1}^2u_{m+1}x_{m+1}^\pm)}\, \alpha^{(m)}(x_{m+1},t),
\\ &&
\beta'(w,s_{m+1}t)=\kappa\frac{\Gamma(t^2s_{m+1}^2,t^2s_{m+1}u_{m+1}w^\pm)}
{\Gamma(s_{m+1}^2,t^2,s_{m+1}u_{m+1}w^\pm)}
\nonumber \\ && \makebox[4em]{} \times
\int_{\mathbb{T}} \frac{\Gamma(s_{m+1}w^\pm x^\pm,u_{m+1}x^\pm)}
{\Gamma(x^{\pm2},t^2s_{m+1}^2u_{m+1}x^\pm)}\,\beta^{(m)}(x,t)\frac{dx}{x}
\nonumber \\ && \makebox[4em]{}
=\kappa\int_{\mathbb{T}}\Gamma(s_{m+1}tw^\pm x_{m+1}^\pm)
\,\alpha'(x_{m+1},s_{m+1}t)\frac{dx_{m+1}}{x_{m+1}}.
\label{ident4}\end{eqnarray}
We substitute into the latter relation expressions (\ref{a_m}) and
(\ref{b_m}) with functions (\ref{a_1}) and (\ref{b_1}) serving
as $\alpha^{(0)}(w,t)$ and $\beta^{(0)}(w,t)$.
As a result, we derive the identity
\begin{eqnarray}\nonumber
\lefteqn{ \int_{\mathbb{T}}\frac{\Gamma(s_{m+1}w^\pm x^\pm,u_{m+1}x^\pm)}
{\Gamma(x^{\pm2},t^2s_{m+1}^2u_{m+1}x^\pm)}
\prod_{k=1}^m\frac{\Gamma(t\prod_{l=k+1}^ms_lu_kx^\pm)}
{\Gamma(t\prod_{l=k+1}^ms_ls_k^2u_kx^\pm)} } &&
\\ &&
\times \frac{\prod_{r=0}^2\Gamma(t\prod_{k=1}^ms_k t_rx^\pm)}
{\Gamma(t\prod_{k=1}^ms_kt_0t_1t_2x^\pm)}\frac{dx}{x}
=\frac{\kappa^m\Gamma(s_{m+1}^2,t^2,s_{m+1}u_{m+1}w^\pm)}
{\Gamma(t^2s_{m+1}^2,t^2\prod_{k=1}^ms_k^2,t^2s_{m+1}u_{m+1}w^\pm)}
\nonumber \\ &&
\times \prod_{k=1}^m\frac{\Gamma(\prod_{l=k}^ms_l^2t^2)}
{\Gamma(s_k^2,\prod_{l=k+1}^ms_l^2t^2)}\prod_{0\leq r<j\leq 2}
\frac{\Gamma(t^2\prod_{k=1}^ms_k^2 t_rt_j)}{\Gamma(t_rt_j)}\,
\nonumber \\ &&
\times \int_{\mathbb{T}^{m+1}}
\frac{\Gamma(ts_{m+1}w^\pm x_{m+1}^\pm,tu_{m+1}x_{m+1}^\pm)
\prod_{r=0}^2\Gamma(t_rx_1^\pm)}{\Gamma(ts_{m+1}^2u_{m+1}x_{m+1}^\pm,
x_1^{\pm 2}, t^2\prod_{k=1}^ms_k^2t_0t_1t_2x_1^\pm)}
\nonumber \\ &&\makebox[1em]{} \times
\prod_{k=1}^m\frac{\Gamma(u_kx_{k+1}^\pm,
\prod_{l=k+1}^ms_l^2t^2s_ku_kx_k^\pm,s_kx_{k+1}^\pm x_k^\pm)}
{\Gamma(x_{k+1}^{\pm2},\prod_{l=k}^ms_l^2t^2u_kx_{k+1}^\pm,s_ku_kx_k^\pm)}
\frac{dx_1}{x_1}\cdots\frac{dx_{m+1}}{x_{m+1}}.
\label{identfin}\end{eqnarray}
For $m=1$ this equality is reduced to relation (\ref{ident1}).

Evidently, a more thorough analysis of consequences of the derived
elliptic Bailey tree for integrals with different multiplicities
of integrations
is desirable. In particular, uniqueness of this tree is an interesting
question for consideration. Implications for $q$-hypergeometric integrals
and elliptic and $q$-hypergeometric series are worth of detailed
investigation as well.

This work is supported in part by the Russian Foundation for Basic
Research (RFBR) Grant No. 03-01-00780.

\end{document}